\documentclass{ifacconf}

\usepackage{graphicx}      % include this line if your document contains figures
\graphicspath{{figures/}}

\usepackage{natbib}        % required for bibliography
\usepackage{comment}        %for commenting out a chunk of content

\usepackage{amsmath} % assumes amsmath package installed
\usepackage{amssymb}  % assumes amsmath package installed
\usepackage{mathtools}

\usepackage{mathrsfs}
\usepackage{xcolor}

\usepackage[width=.75\textwidth]{caption}
\usepackage{tabularx}
\usepackage{booktabs}           % More beautiful tables

\usepackage[per-mode=symbol]{siunitx}

\usepackage{pgfplots}\pgfplotsset{compat=1.9}
\usepackage{grffile}
\pgfplotsset{compat=newest}  % Use advanced features
\usetikzlibrary{plotmarks}
\usetikzlibrary{arrows.meta}
\usepgfplotslibrary{patchplots}

\usepackage{subcaption}

\usepackage{todonotes}
\newcommand{\todoVA}[1]{\todo[inline,color=yellow!60!white]{VA: #1}}

\begin{document}
\begin{frontmatter}

\title{Safety-Aware and Data-Driven Predictive Control for Connected Automated Vehicles at a Mixed Traffic Signalized Intersection\thanksref{footnoteinfo}}

\thanks[footnoteinfo]{This research was supported by ARPAE's NEXTCAR program under the award number DE-AR0000796.}

\author[First]{A M Ishtiaque Mahbub},
\author[First]{Viet-Anh Le},
\author[First]{Andreas A. Malikopoulos}

\address[First]{University of Delaware,
   Newark, DE 19716 USA (e-mail: mahbub@udel.edu, vietale@udel.edu, andreas@udel.edu).}
%\address[Second]{Colorado State University,
%    Fort Collins, CO 80523 USA (e-mail: author@lamar. colostate.edu)}
% \address[Third]{Electrical Engineering Department,
%    Seoul National University, Seoul, Korea, (e-mail: author@snu.ac.kr)}

\begin{abstract}                % Abstract of not more than 250 words.
A typical urban signalized intersection poses significant modeling and control challenges in a mixed traffic environment consisting of connected automated vehicles (CAVs) and human-driven vehicles (HDVs). In this paper, we address the problem of deriving safe trajectories for CAVs in a mixed traffic environment that prioritizes rear-end collision avoidance when the preceding HDVs approach the yellow and red signal phases of the intersection. We present a predictive control framework that employs a recursive least squares algorithm to approximate in real time the driving behavior of the preceding HDVs and then uses this approximation to derive safety-aware trajectory in a finite horizon. We validate the effectiveness of our proposed framework through numerical simulation and analyze the robustness of the control framework.
\end{abstract}


\end{frontmatter}
%===============================================================================

\section{Introduction}
 
Optimal coordination of connected automated vehicles (CAVs) can improve network performance, e.g., fuel economy, traffic throughput, at traffic scenarios such as urban intersections, see \cite{talebpour2016CAVTraffic}.
In recent efforts, a decentralized optimal control framework has been established for real-time coordination of CAVs traveling through signal-free automated intersections, see \cite{Malikopoulos2020,chalaki2020TITS, mahbub2020Automatica-2, Kumaravel:2021uk}. 
%
%mixed traffic
However, these approaches have been developed based on the strict assumption of 100\% CAV penetration rate which is currently not realizable;  see \cite{alessandrini2015automatedmixed2060}. 

CAVs must be able to safely co-exist with human-driven vehicles (HDVs) resulting in a \emph{mixed traffic environment}, which pose significant modeling and control challenges due to the stochastic nature of human-driving behavior. 
Furthermore, the use of conventional traffic lights is still the most prevalent way of traffic control at urban intersections that adds an extra layer of complexity in modeling the HDV behavior due to the presence of the zone of yellow light dilemma; see \cite{zhang2014yellowdilemma}. 
Thus, the need for an efficient CAV control framework considering the inclusion and interaction of HDVs approaching the signalized intersections is essential to provide safety assurance under unknown HDV behavior.

Several research efforts have adopted adaptive cruise control (ACC) for automated vehicles in a mixed traffic environment to tackle the HDV behavior and ensure rear-end collision avoidance; see \cite{jiang2007mixedCACC, yuan2009trafficMixedACC}. 
\cite{lu2019ecological} considered a variation of the car-following model to design an eco ACC controller.
%%ACC/CACC
However, ACC controllers using car-following models such as the intelligent driver model (IDM) \citep{treiber2013traffic} do not always perform well since they can have stability implications leading to rear-end collision; see \cite{milanes2014ACC}.
\cite{milanes2013cooperative} proposed a cooperative adaptive cruise controller where the control parameters are derived using system identification on real-world experimental data. However, such control parameters cannot capture the instantaneous changes in HDV behavior.
\cite{naus2010mpcAcc} proposed an explicit model predictive control (MPC) ACC controller that employs a prediction model with a constant speed assumption of the preceding vehicle and does not consider the complex car-following dynamics of the human driver.
\cite{dollar2021mpc} utilized an IDM model to identify offline the human driving styles in a car-following scenario and developed an MPC-based cruise control for a CAV.
\cite{jin2018connected} proposed an optimal cruise control design in which feedback gains and driver reaction time of HDVs were estimated in real time by a sweeping least square method.

%%%%%%%%%%%%%%%%%%%%%%%%%%%%%%%%%%%%%%%%%%%%%%%%%%%%%%%%%%%%%%%%
%contribution of the paper
In this paper, we consider the problem related to controlling a CAV while approaching a signalized intersection in the presence of multiple preceding HDVs with unknown driving behavior. 
To generate safe and optimal control actions for the CAV, we propose a data-driven predictive control framework that takes into account the future trajectories of the HDVs to ensure that the collision does not take place over a finite-time horizon.
The control framework is then implemented in a receding horizon manner for robustness against stochastic driving behavior of HDVs.
The constant time headway relative velocity (CTH-RV) model and the recursive least squares (RLS) algorithm are utilized to estimate the HDV's driving behavior given the data collected online.
We evaluate the efficiency of the proposed method by numerical simulations that employ a nonlinear car-following model to replicate the human drivers.

%%%%%%%%%%%%%%%%%%%%%%%%%%%%%%%%%%%%%%%%%%%%%%%%%%%%%%%%%%%
%Structure of the paper
The remainder of the paper is organized as follows. 
In Section \ref{sec:pf}, we present the modeling framework and formulate the problem. 
In Section \ref{sec:control_framework}, we provide a detailed exposition of the safety-aware and data-driven predictive control framework with real-time behavior estimation. 
In Section \ref{sec:result}, we evaluate the effectiveness of the proposed approach in a simulation environment. 
Finally, we draw conclusions and discuss the future research directions in Section \ref{sec:conclusion}.
%
%%%%%%%%%%%%%%%%%%%%%%%%%%%%%%%%%%%%%%%%%%%%%%%%%%%%%%%%%%%%%%%%%
%
%%%%%%%%%%%%%%%%%%%%%%%%%%%%%%%%%%%%%%%%%%%%%%%%%%%%%%%%%%%%%%%%%
\section{Problem Formulation}\label{sec:pf}
%!TEX root = main.tex

\begin{figure}
\begin{center}
\includegraphics[scale = 0.7,trim={6.5cm 8cm 7cm 7.5cm},clip]{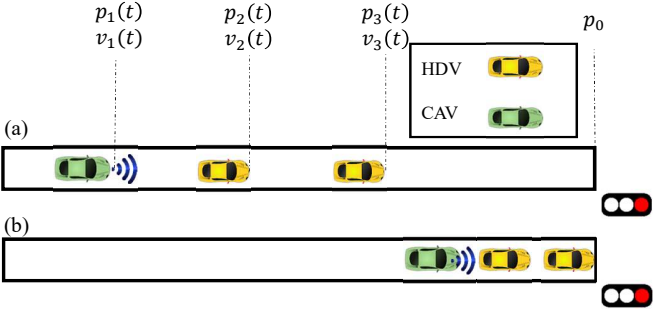}    % The printed column width is 8.4 cm.
\caption{A connected and automated vehicle (green) and human-driven vehicles (yellow) approaching a signalized intersection at red signal phase.}
\label{fig:problem_formulation}
\end{center}
\end{figure}

We consider multiple HDVs followed by a CAV traveling on a single-lane road and approaching an urban signalized intersection with a red (or yellow then red) traffic signal phase (Fig. \ref{fig:problem_formulation}).
Note that, the general idea of our formulation can be extended to different cases such as yield/stop traffic sign, downstream traffic congestion, and pedestrian crossing, where the preceding HDVs' motion can change abruptly to come to a full stop.
Next, to facilitate our exposition we provide the following definitions.
\begin{defn}\label{def:N(t)}
Suppose that the red signal phase is active at some time instant $t = t^0$. The set $\mathcal{N}$ of the vehicles approaching the intersection at $t = t^0$ is $\mathcal{N}=\{N,N-1, \ldots, 1\}$, where $N\in \mathbb{N}$ is the total number of vehicles under consideration. Here, the vehicles are assigned integer indices by the order of their respective distances from some fixed stopping position $p_0$ located downstream near the signal head. The indices $N, N-1, \ldots, 2$ represent the HDVs followed by the CAV denoted by the index $1$.
\end{defn}
\begin{defn}\label{def:N_hdv}
The set of HDVs at time instant $t = t^0$ is $\mathcal{N}_{\text{HDV}}=\mathcal{N} \setminus \{1\}$.
\end{defn}

When the red signal phase is active at some time instant $t = t^0$, the HDV-$N$ in $\mathcal{N}$ must stop behind the position $p_0$. The objective of the CAV-$1$ is to derive an optimal trajectory to stop behind HDV-$2$ such that no rear-end collision takes place.
\begin{rem}\label{rem:empty_N_HDV}
In our formulation, we require that the set of HDVs $\mathcal{N}_{\text{HDV}}$ is non-empty at time instant $t = t^0$ when the red signal phase of the intersection is active. If $\mathcal{N}_{\text{HDV}}$ is empty, then the problem of avoiding rear-end collision becomes redundant.
\end{rem}
%
%Communication
\subsection{Communication Topology}\label{subsec:communication_topology}
The CAV-$1$ is retrofitted with appropriate sensors and communication devices to estimate in real time the state information of the preceding HDVs in $\mathcal{N}_{\text{HDV}}$ through vehicle-to-everything communication protocol and intelligent roadside units; see \cite{duan2021_VideoAsSensor}. We refer to the unidirectional flow of information from the preceding HDVs in $\mathcal{N}_{\text{HDV}}$ to the trailing CAV-$1$ as the multi-predecessor communication topology. 

We impose the following assumption.
\begin{assum}	\label{assum:communication}
Communication to and from the CAV occurs without any delays and errors.
\end{assum}
Assumption \ref{assum:communication} may be strong, but it is relatively straightforward to relax it as long as the noise in the measurements and/or delays is bounded. %For example, we can determine upper bounds on the state uncertainties as a result of sensing or communication errors and delays, and incorporate these into more conservative safety constraints.

In contrast, the trajectory of each HDV $i$ in $\mathcal{N}_{\text{HDV}}$ is solely dictated by the perception of the state information of the immediate preceding HDV $i+1$ in $\mathcal{N}_{\text{HDV}}$.
For the leading vehicle HDV-$N$ that does not have a preceding vehicle, its driving actions depend on the relative distance to the stopping point.

\subsection{Vehicle Dynamics and Constraints}
We consider the following discrete-time model with a sampling time $\tau\in\mathbb{R}^+$ to represent the dynamics of each vehicle $i\in\mathcal{N}$,
\begin{subequations}\label{eq:dynamics_pv}
\begin{align}
    &p_i(t+1) = p_i(t) + v_i(t) \tau + \frac{1}{2}u_i(t) \tau^2,\\
    &v_i(t+1) = v_i(t) + u_i(t) \tau,
\end{align}
\end{subequations}
where $p_{i}(t)\in\mathcal{P}_{i}$, $v_{i}(t)\in\mathcal{V}_{i}$, and $u_{i}(t)\in\mathcal{U}_{i}$ denote the position, speed and control input (acceleration/deceleration) of each vehicle $i$ in $\mathcal{N}$.
The sets $\mathcal{P}_{i}$, $\mathcal{V}_{i}$, and $\mathcal{U}_{i}$, $i\in\mathcal{N}(t),$ are complete and totally bounded subsets of $\mathbb{R}$.
Note that in the discrete-time dynamics model \eqref{eq:dynamics_pv}, we assume that the control input $u_i(t)$ of each vehicle $i$ in $\mathcal{N}$ remains constant in the time period of length $\tau$ between time instants $t$ and $t+1$, which is different to some previous approaches that assume constant speed between time instants $t$ and $t+1$; see \cite{naus2010mpcAcc, kianfar2012MPC_ACCChallenge}.

%constraints
To ensure that the control input and vehicle speed are within a given admissible range, the following constraints are imposed,
\begin{subequations}\label{eq:speed_accel constraints}
\begin{align}
u_{\mathrm{min}}   &\leq u_{i}(t)\leq u_{\mathrm{max}},\quad\text{and}\\
0   \leq v_{\mathrm{min}}&\leq v_{i}(t)\leq v_{\mathrm{max}},
\end{align}
\end{subequations}
where $u_{\mathrm{min}}$, $u_{\mathrm{max}}$ are the maximum braking and acceleration, respectively, of each vehicle in $\mathcal{N}$, and $v_{\mathrm{min}}$, $v_{\mathrm{max}}$ are the minimum and maximum speed limits, respectively.

%control structure u_i(t)
%car-following dynamics
The control input $u_i(t)$ of each vehicle $i\in\mathcal{N}$ in \eqref{eq:dynamics_pv} can take different forms based on the consideration of connectivity and automation.
%CAV dynamics
For CAV-$1$ in $\mathcal{N}$, we consider a switching control framework based on the following cases: if at time instant $t=t^0$ (Remark \ref{rem:empty_N_HDV}) (a) $\mathcal{N}_{\text{HDV}}$ is empty, then CAV-$1$ derives its control input by using its default adaptive cruise controller (see \cite{milanes2014ACC}), (b) if $\mathcal{N}_{\text{HDV}}$ is not empty, then CAV-$1$ derives and implements the control input $u_1(t)$ using the proposed control framework discussed in Section \ref{sec:control_framework}.

%HDV dynamics
%
For each HDV $i\in \mathcal{N}_{\text{HDV}}$, however, we consider a car-following model to represent the predecessor-follower coupled dynamics (Fig. \ref{fig:problem_formulation}) with its preceding vehicle $i+1$ that has the following generic structure
\begin{gather}\label{eq:hdv_dynamics}
     {{u}_i(t) = f_i (\Delta p_{i}(t), v_i(t), \Delta v_i(t)),}
\end{gather}
where $f_i(\cdot)$ represents the behavioral function of the car-following model of vehicle $i\in\mathcal{N}_{\text{HDV}}$, and $\Delta p_i(t):=p_{i+1}(t)-p_i(t)-l_c$ and $\Delta v_i(t):=v_{i+1}(t)-v_i(t)$ denote the headway and approach rate of vehicle $i$ with respect to its preceding vehicle $i+1$, respectively.
%edge cases
We consider two edge cases that may arise from the above definitions: (a) if there is no vehicle $i+1$ preceding vehicle $i$ within a certain look-ahead distance $d_f$, then we consider $\Delta p_i(t) = d_f$ and {$\Delta v_i(t)=0$}, and (b) if there is an obstruction/red signal phase immediately ahead of vehicle $i$ at a distance $d_s$, then $\Delta p_i(t) = d_s$ and $\Delta v_{i}(t)=-v_i(t)$.
%CFM Literature
There are several car-following models reported in the literature that can emulate a varied class of human driving behavior; see \cite{CFM_VehicularTraffic}.
% Given the dynamics model in \eqref{eq:dynamics_pv} and the car-following model \eqref{eq:hdv_dynamics}, we can predict the future states of the HDV and use the prediction in the control framework.

% VA: The following sentence my conflict with the fact that we use CFM to predict the state of the HDV, so I changed it
% In our formulation, the CAV does not have any knowledge of the behavioral function $f(\cdot)$ of the preceding HDVs.
% \todoVA{I think the following paragraph should be moved to Section 3, in the RLS subsection}
The parameters of a car-following model can be recovered from historical data using offline identification methods; see \cite{treiber2013traffic}.
However, since the historical data might not be available and the human driving behavior usually changes over time, offline identification methods do not work well in practice.
As a result, in our proposed framework, we consider that the CAV does not have full prior knowledge of the behavioral function $f_i(\cdot)$ of the preceding HDVs.
Instead, the CAV assumes a specific type of car-following model for the HDV, then estimates the model parameters for each HDV online using real-time collected data.
A method for estimating car-following model parameters of the HDVs is given in Section~\ref{sec:parameter_estimation}.

%
%information set
%tracking error states
To capture the car-following characteristics of the preceding HDV-$2$'s dynamics from the CAV-$1$'s control point of view, we define additional states as %associated with tracking the headway $e_{p}(t)$ and speed deviation ${e}_{v}(t)$ as
 \begin{subequations}\label{eq:tracking_error}
 \begin{align}
    &e_{p}(t) = p_2(t)-p_1(t)-l_c,\\
    &e_{v}(t) = v_2(t)-v_1(t).
\end{align}
\end{subequations}
% In this paper, we consider that the control input $u_i(t)$ of the HDVs in \eqref{eq:hdv_dynamics} remain constant between the time instances $k$ and $k+1$. Therefore, we can use the dynamics model in \eqref{eq:dynamics_pv} to predict the future states,

%rear-end collision
To introduce the rear-end collision avoidance constraint, we first use the following definition of dynamic safe following headway $s_i(t)$.
\begin{defn}\label{defn:s_i}
The dynamic safe following headway $s_i(t)$ between two consecutive vehicles $i \text{ and }(i+1)\in\mathcal{N}$ is
\begin{equation}\label{eq:s_i}
    { s_i(t)= \rho_iv_i(t)+ s_0,}
\end{equation}
where $\rho_i\in \mathbb{R}^+$ denotes a desired time headway that each vehicle $i\in\mathcal{N}$ maintains while following the preceding vehicle, and $s_0\in \mathbb{R}^+$ is the standstill distance denoting the minimum bumper-to-bumper gap at stop.
\end{defn}
The rear-end collision avoidance constraint between CAV-$1$ and its immediately preceding HDV-$2$ can thus be written as
\begin{equation}
e_p(t) \ge s_1(t).
\label{eq:rearend_constraint}
\end{equation}
%
%Let $x_{i}(t)=\left[p_{i}(t) ~ v_{i}(t)\right]  ^{T}$ denote the state of each vehicle $i$, with initial value $x_{i}^{0}=x_i(t_i^0)=\left[p_{i}^{0} ~ v_{i}^{0}\right]  ^{T}$, where $p_{i}^{0}= p_{i}(t_{i}^{0})=0$ at the entry of the corridor, taking values in $\mathcal{X}_{i}=\mathcal{P}_{i}\times\mathcal{V}_{i}$. The state space $\mathcal{X}_{i}$ for each vehicle $i$ is closed with respect to the induced topology on $\mathcal{P}_{i}\times \mathcal{V}_{i}$ and thus, it is compact.
%

\begin{comment}
\todoVA{Is the following assumption important? If not, we can omit to save more space.}
%assumption
In our analysis, we impose the following assumption:
\begin{assum}\label{assum:feasible}
At $k = t^0$, none of the state, control and safety constraints are active for CAV $1$.
\end{assum}
%
Assumption \ref{assum:feasible} ensures that the CAV-$1$ derives the optimal control input starting from a feasible state space.
\end{comment}

We now formalize the main objective of the CAV-$1$ control framework.
\begin{prob}\label{problem:objectective}
{Given the multi-predecessor communication topology (Section \ref{subsec:communication_topology}), the main objective of CAV-$1$ is to derive its optimal control input $u_1^*(t)$ such that CAV-$1$ adapts to its preceding HDV's driving behavior in real time and drives the states $e_p(t)$ and $e_v(t)$ to their respective reference states with minimum control effort satisfying the state, control and safety constraints in \eqref{eq:speed_accel constraints}-\eqref{eq:rearend_constraint}.}
%estimate the car-following parameters in \eqref{eq:hdv_dynamics} online, and
\end{prob}
%%%%%%%%%%%%%%%%%%%%%%%%%%%%%%%%%%%%%%%%%%%%%%%%%%%%%%%%%%%%%%%%%
%
%%%%%%%%%%%%%%%%%%%%%%%%%%%%%%%%%%%%%%%%%%%%%%%%%%%%%%%%%%%%%%%%%
\section{Control Framework} \label{sec:control_framework}
%!TEX root = main.tex

% \todoVA{I would suggest to use "predictive control" throughout the paper because s main feature is that we predict the HDV's future trajectories.}

\begin{figure}
    \centering
    \includegraphics[scale = 0.6,trim={6cm 7cm 7cm 5cm},clip]{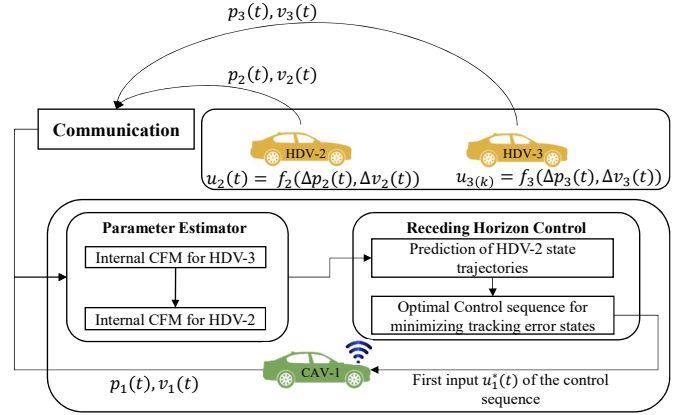}
    \caption{The structure of the proposed control framework to address Problem \ref{problem:objectective}. }
    \label{fig:control_framework}
\end{figure}

In our approach, we adopt a receding horizon predictive control framework with multi-predecessor communication topology and data-driven estimation of HDVs' car-following parameters for state prediction to address Problem \ref{problem:objectective}, as shown in Figure \ref{fig:control_framework}. %MPC
In the receding horizon control, the optimal control input at the current time step is obtained by solving a predictive control problem with a horizon $T_p$ while only the first element of the obtained control input sequence is implemented. 
Afterward, the horizon moves forward one step, and the above process is repeated until a final horizon is reached; see \cite{borrelli2017predictive}.
{Note that, the prediction horizon $T_p$ is usually selected empirically to best accommodate the control performance and computational requirement.} The essential steps of the proposed framework are outlined as follows.
%trim={<left> <lower> <right> <upper>}
%
\begin{enumerate}
    \item \textbf{Data-driven parameter estimation:} At each time instant $t$, the current states $p_i(t), v_i(t)$ of each preceding HDV $i$ in $\mathcal{N}_{\text{HDV}}$ is communicated to CAV-$1$. Since the exact car-following model $f_i$ of each HDV $i$ in $\mathcal{N}_{\text{HDV}}$ is unknown to CAV-$1$, it considers a specific type of car-following model to represent the driving behavior of each HDV, and estimates the parameters of the car-following model for each HDV online. 
    \item \textbf{Predictive control problem:} CAV-$1$ then uses the estimated car-following model from Step 1 to predict the future state trajectories of the immediately preceding HDV-$2$ and derives its own optimal control input sequence ${U}_1^*(t):=[{u}_1^*(t),~{u}_1^*(t+1),\ldots,{u}_1^*(t+T_p-1)]^T$ using the receding horizon control framework discussed above. Finally, CAV-$1$ implements only the first control input $u_1^*(t)$.
\end{enumerate}
In what follows, we provide a detailed exposition of the steps discussed above.
\subsection{Online Car-following Model Parameter Estimation}\label{sec:parameter_estimation}
In this section, we use a recursive least-squared formulation \citep{ljung1983RLS} to estimate the parameters of the internal car-following model residing in CAV-$1$'s mainframe to represent the driving behavior of each of the preceding HDVs. To this end, we consider the CTH-RV model \citep{wang2020onlineRLS}
\begin{equation}\label{eq:cth-rv}
\begin{multlined}
    v_i(t+1) = v_i(t) + \eta_i(\Delta p_i(t) - \rho_i v_i(t))\tau + \\
    \nu_i (v_1(t)-v_i(t))\tau,
\end{multlined}
\end{equation}
where the model parameters $\eta_i$ and $\nu_i$ are the control gains on the constant time headway and the approach rate, and $\rho_i$ is the desired safe time headway for each HDV $i$ in $\mathcal{N}_{\text{HDV}}$, respectively.
We employ the linear CTH-RV model instead of other complex nonlinear models so that the resulting control problem presented in the next section is thus convex and can be solved efficiently in real-time.
Moreover, it is also observed that the CTH-RV model is highly comparable to other nonlinear car-following models in terms of data fitting \citep{gunter2019modeling}. 

%\todoVA{Why do we need to have $\tau_e$ here? Why don't we keep it just $\tau$?}

Suppose that we measure the speed $v_i(t)$, headway gap $\Delta p_i(t)$ and approach rate $\Delta v_i(t)$ for each preceding HDV $i$ in $\mathcal{N}_{\text{HDV}}$ with sampling rate $\tau$. We recast \eqref{eq:cth-rv} as
\begin{align}\label{eq:estimation_1}
    v_i(t+1) = \gamma_{i,1} v_i(t) + \gamma_{i,2} \Delta p_i(t) + \gamma_3 v_1(t),
\end{align}
where $\gamma_{i,1}:= (1-(\eta_i \rho_i+\nu_i)\tau)$, $\gamma_{i,2}:= \eta_i \tau$ and $\gamma_{i,3}:= \nu_i \tau$
are the parameters that can be estimated using the RLS algorithm.
The original model parameters $\eta_i, \nu_i$ and $\rho_i$ are then uniquely determined from $\gamma_{i,1}, \gamma_{i,2}, \gamma_{i,3}$ as long as $\gamma_{i,2} \neq 0$.
Next, we can write \eqref{eq:estimation_1} in matrix form as
\begin{align}
    v_i(t+1) = \gamma_i^T \phi_i(t),
\end{align}
where $\phi_i(t):=[v_i(t),~ \Delta p_i(t), ~v_1(t)]^T$ is the regressor vector and $\gamma_i:=[\gamma_{i,1},~\gamma_{i,2},~\gamma_{i,3}]^T$ is the parameter vector. 
We can estimate $\gamma_i$ using the following recursive least squares algorithm as follows \citep{ljung1983RLS}
\begin{subequations}
\begin{align}\label{eq:rls}
    &\hat{\gamma_i}(t) = \hat{\gamma_i}(t-1)+ L_i(t)[v_i(t)-\hat{v}_i(t)], \\
    & \hat{v}_i(t) = \hat{\gamma_i}^T(t-1)\phi_i(t), \\
    & L_i(t) = \frac{P_i(t-1)\phi_i(t)}{\xi + \phi_i^T(t)P_i(t-1)\phi_i(t)},\\
    & P_i(t) = \frac{1}{\xi} \bigg[P_i(t-1) - \frac{P_i(t-1)\phi_i(t)\phi_i^T(t)P_i(t-1)}{\xi+ \phi_i^T(t)P_i(t-1)\phi_i(t)}\bigg ].
\end{align}
\end{subequations}
Here, $\xi\in[0,1]$ is the forgetting factor that assigns a higher weight to the recently collected data points and discounts older measurements, and $\hat{\gamma}_i (t)$ denotes the estimate of the parameter vector $\gamma_i$ at time instant $t$, which is updated recursively as new data becomes available.
%\todoVA{We should explain more about the terms in RLS.}
%
In what follows, we introduce the predictive control problem that is needed to be solved. % in a receding horizon manner.
\subsection{Predictive Control Problem}\label{subsec:mpc}
The main objective of the predictive controller of the CAV is to (a) drive the position tracking state $e_p(t)$ to a reference $e_{p,r}(t)$, (b) drive the speed tracking state $e_v(t)$ to zero, and (c) minimize CAV-$1$'s control input $u_1(t)$. 
To this end, the receding horizon controller generates the predictive states $e_p(t+n|t), e_v(t+n|t)$ for $n=1,\ldots, T_p$ at each time instant $t$ for a predictive horizon $T_p$ using the state definitions in \eqref{eq:tracking_error}, vehicle dynamics in \eqref{eq:dynamics_pv} and internal car-following models of the HDVs in \eqref{eq:cth-rv} approximated in the previous section. 
Then the control input sequence ${U}_1(t):=[{u}_1(t),~{u}_1(t+1),\ldots,{u}_1(t+T_p-1)]^T$ is derived such that the predictive states are driven to their respective reference states. The predictive control problem thus can be written as
\begin{align}\label{eq:ocp}
     \min_{U_1(t)} \frac{1}{2} &\sum_{n=1}^{T_p} \bigg[ w_{e_p}(e_p(t+n|t) - e_{p,r}(t+n|t))^2 \\&+ w_{e_v}e_v(t+n|t)^2 +  w_{u} (u_1(t+n-1))^2\bigg],\nonumber\\
&\text{subject to}:\nonumber\\
&\text{model: } \eqref{eq:dynamics_pv}, \eqref{eq:hdv_dynamics}, \eqref{eq:tracking_error},\nonumber\\
&\text{constraints: } \eqref{eq:speed_accel constraints}, \eqref{eq:rearend_constraint}, \nonumber\\
&\text{reference state: } e_{p,r}(t):=s_1(t), \nonumber
\end{align}
where the predictive reference state $e_{p,r}(t+n|t)$ can be computed using the relation $e_{p,r}(t)=s_1(t)$ and the dynamics model in \eqref{eq:dynamics_pv} and \eqref{eq:tracking_error}, and $w_{e_p}, w_{e_v}, w_u \in \mathbb{R}^+$ are the weights on the reference tracking of the headway $e_p(t)$, speed deviation $e_v(t)$, and the CAV-$1$'s control input $u_1(t)$, respectively.

The predictive control problem in \eqref{eq:ocp} can be transformed into a standard constrained quadratic programming problem and solved using commercially available solvers; see \cite{CasAdi_python}. At each discrete time instant $t$, the optimal control sequence $U_1^*(t)$ is computed by solving \eqref{eq:ocp} and only the first control input $u_1^*(t)$ is applied. Then the system moves to the next time instant $t+1$ and the process is repeated until a final time horizon is reached.

\begin{rem}\label{rem:updating_vehicleSet}
While implementing the above control framework, if any of the preceding HDVs leaves the current lane or passes the intersection at any time instant $t$, we simply update the sets $\mathcal{N}$ and $\mathcal{N}_{\text{HDV}}$ starting from the next time instant $t+1$, where the control problem \eqref{eq:ocp} is again solved with the updated information.
\end{rem}
%%%%%%%%%%%%%%%%%%%%%%%%%%%%%%%%%%%%%%%%%%%%%%%%%%%%%%%%%%%%%%%%%%%%%%%%%%%%%%%%%%%%%%
%
%%%%%%%%%%%%%%%%%%%%%%%%%%%%%%%%%%%%%%%%%%%%%%%%%%%%%%%%%%%%%%%%%%%%%%%%%%%%%%%%%%%%%%
\section{Simulation and Result}\label{sec:result}
%!TEX root = main.tex

This section validates the performance of the proposed safety-aware data-driven predictive control by numerical simulations at a mixed-traffic signalized intersection.

\subsection{Simulation Setup} 

In the simulations, we utilize a nonlinear car-following model namely the optimal velocity model (OVM) to generate the driving actions of simulated human drivers \citep{bando1995dynamical}.
The car-following OVM is given as
\begin{equation}
\begin{split}
u_i (t) &= \alpha \big( V_i(t) - v_i(t) \big) + \beta \Delta v_i(t), \\
V_i (t) &= \frac{v_{\text{d}}}{2} \big( \tanh \big( \Delta p_i(t)-s_i(t)) + \tanh(s_i(t) \big) \big).
\end{split}
\end{equation}
The parameters of the OVM for each HDV include the driver’s sensitivity coefficients $\alpha$ and $\beta$, and the desired speed $v_d$.
These parameters for the simulated HDVs are assumed to be different to each other and chosen by random perturbations up to 20\% around the following nominal values: $\alpha = 0.8$, $\beta = 0.6$, $v_d = \SI{15.0}{m/s}$, $\rho = \SI{2.0}{s}$, $s_0 = \SI{5.0}{m} $. 
The parameters and weights in the predictive control framework used for the simulations are given in Table~\ref{tab:sim-params}.
The RLS-based estimators are initialized with the following values: 
% $\eta = 1.0$, $\nu = 1.8$, $\tau_e = 1.5$. 
$\gamma_i(0) = [0.67, 0.1, 0.18]^T$ and $P_i(0) = 0.01 \, \mathbb{I}_3$ where $\mathbb{I}_3$ is the $3 \times 3$ identity matrix, while the forgetting factor is chosen as $\xi = 1.0$. The impact of $\xi$ on RLS algorithm is investigated in detail by \cite{vahidi2005recursive} and thus, omitted here.
% \todoVA{An initial estimate of the parameters, $\gamma_i(0)$ and $P_i(0)$}
%
Python is used in the simulations in which the constrained optimal control problem is formulated by CasADi framework; see \cite{andersson2019casadi}, and solved by the built-in qpOASES solver.
%
% Note that once the traffic signal turns to yellow, the coordinator starts collecting data from the HDVs and then transmit the data to the CAV while the CAV starts estimating the car-following parameters of the preceding HDVs and using the satey-aware predictive controller.
%
\begin{table}[!bt]
  \caption{Parameters of the controller}
  \label{tab:sim-params} 
  \centering
  \begin{tabular}{p{0.085\textwidth}p{0.09\textwidth}p{0.085\textwidth}p{0.09\textwidth}}
    \toprule[1pt]% <-- Toprule here
    \textbf{Parameters} & \textbf{Value} & \textbf{Parameters} & \textbf{Value} \\
    \midrule[0.5pt] % <-- Midrule here
    $\tau$ & \SI{0.1}{s} & $T_p$ & 50 \\
    $v_{\mathrm{max}}$ & \SI{15}{m/s} & $v_{\mathrm{min}}$ & \SI{0}{m/s}\\
    $u_{\mathrm{max}}$ & \SI{3}{m/s^2} & $u_{\mathrm{min}}$ & \SI{-5}{m/s^2}\\
    $\rho$ & \SI{2.0}{s} & $s_0$ & \SI{3.0}{m} \\
    $w_{e_p}$ & 1 & $w_{e_v}$ & 0.1\\
    $w_u$ & 1 & &\\
    \bottomrule[1pt] % <-- Bottomrule here
  \end{tabular}
\end{table}

\subsection{Results and Discussions} 

\begin{figure*}[tb]
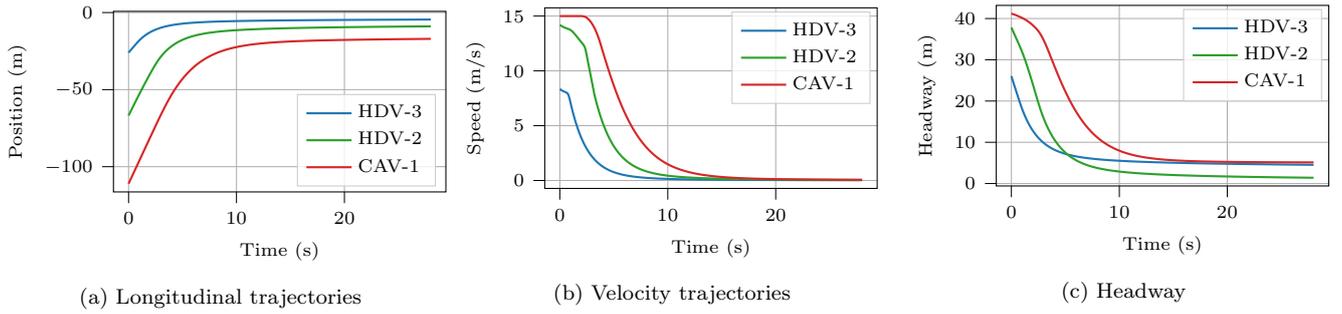

\centering
\begin{subfigure}{.32\textwidth}
    \centering
    \scalebox{0.99}{\input{figures/sim_position.tex}}
    \caption{Longitudinal trajectories}
    \label{fig:sim1_position}
\end{subfigure}
\begin{subfigure}{.32\textwidth}
    \centering
    \scalebox{0.99}{\input{figures/sim_speed.tex}}
    \caption{Velocity trajectories}
    \label{fig:sim1_velocity}
\end{subfigure}
\begin{subfigure}{.32\textwidth}
    \centering
    \scalebox{0.99}{\input{figures/sim_headway.tex}}
    \caption{Headway}
    \label{fig:sim1_headway}
\end{subfigure}
\caption{Longitudinal trajectories, velocities and distances of the vehicles in the simulation with 3 vehicles.}
\label{fig:sim1}
\end{figure*}

\begin{figure*}[tb]
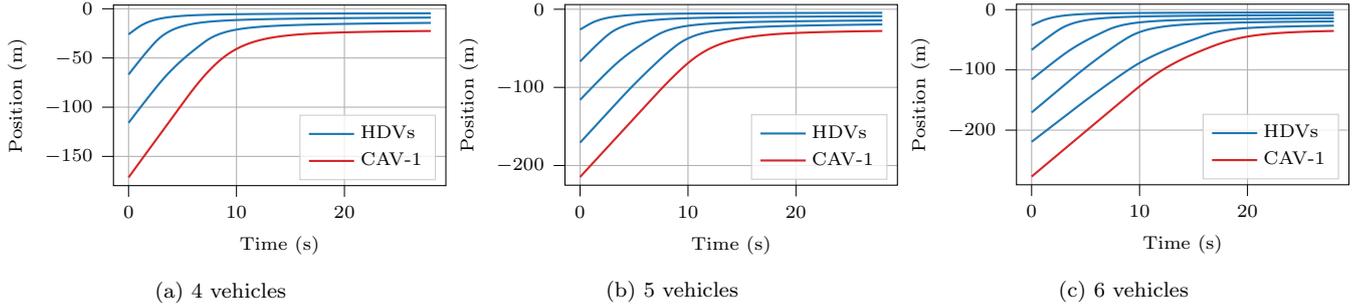

\centering
\begin{subfigure}{.32\textwidth}
    \centering
    \scalebox{0.99}{\input{figures/sim_position4.tex}}
    \caption{4 vehicles}
    \label{fig:sim4}
\end{subfigure}
\begin{subfigure}{.32\textwidth}
    \centering
    \scalebox{0.99}{\input{figures/sim_position5.tex}}
    \caption{5 vehicles}
    \label{fig:sim5}
\end{subfigure}
\begin{subfigure}{.32\textwidth}
    \centering
    \scalebox{0.99}{\input{figures/sim_position6.tex}}
    \caption{6 vehicles}
    \label{fig:sim6}
\end{subfigure}
\caption{Longitudinal trajectories in the simulations with different numbers of vehicles approaching the intersection.}
\label{fig:sim_multiple}
\end{figure*}
The results for a numerical simulation involving a CAV and 2 preceding HDVs are illustrated in Fig.~\ref{fig:sim1}, in which the longitudinal positions, speeds, and headways of all the vehicles are given in Figures~\ref{fig:sim1_position}, \ref{fig:sim1_velocity} and \ref{fig:sim1_headway}, respectively.
Note that the position by which the vehicles must stop is $p_0 = 0$, and for the leading HDV, the headway is computed as the relative distance to the stopping point. 
As can be seen from Figures~\ref{fig:sim1_position}-\ref{fig:sim1_headway}, the simulated HDVs slow down and then stop while approaching the signalized intersection.
Given the behavior of the HDVs, the proposed control framework can perform safe and comfortable braking for the CAV without violating any of the state, input, and safety constraints.    

Moreover, to assess the scalability of the proposed control framework to the number of preceding vehicles, we conduct three other simulations for the scenarios with 4, 5, and 6 vehicles (3, 4, and 5 HDVs, respectively) and illustrate the vehicle trajectories in Fig.~\ref{fig:sim_multiple}. 
These results verify that the proposed control framework works effectively with different numbers of preceding vehicles. 
% \revise{More simulation results can be found in \url{https://sites.google.com/view/ud-ids-lab/DPC-Signalized-Intersection}.}

\begin{figure}[!tb]
    \centering
    \scalebox{0.95}{\input{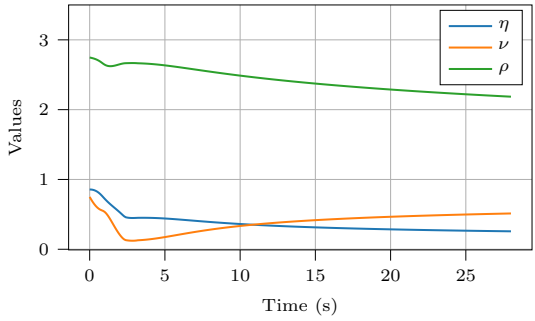}}
    \caption{Estimates for the car-following parameters of the HDV-2}
    \label{fig:estimation}
    \vspace{-10pt}
\end{figure}

Finally, the estimated parameters in the CTH-RV car-following model for the HDV-2 are depicted in Fig.~\ref{fig:estimation}.
{As more real-time data are added to update the estimations, the car-following parameters stabilize to the set of values that accurately describes the driving behavior of the HDVs.}
Therefore, using the linear CTH-RV model and online RLS technique, we can approximate a nonlinear car-following model and use this estimation to predict the future states of the HDVs.   
% The convergence rate of the RLS algorithm depends on \textcolor{red}{something}.

% \todoVA{I am thinking the convergence rate of the RLS depends on what.}

%%%%%%%%%%%%%%%%%%%%%%%%%%%%%%%%%%%%%%%%%%%%%%%%%%%%%%%%%%%%%%%%%%%%%%%%%%%%%%%%%%%%%%%%%%%
%
%%%%%%%%%%%%%%%%%%%%%%%%%%%%%%%%%%%%%%%%%%%%%%%%%%%%%%%%%%%%%%%%%%%%%%%%%%%%%%%%%%%%%%%%%%%%
%
%
\section{Concluding Remarks}\label{sec:conclusion}
In this paper, we addressed the problem of a CAV traveling in a mixed traffic environment and approaching a signalized intersection.
A data-driven predictive control framework was developed in which the car-following behavior of HDVs ahead of the CAV is modeled by the CTH-RV model with online estimated parameters through the RLS algorithm.
In the proposed framework, by utilizing data-driven car-following models, the CAVs can predict the future behavior of the HDVs and then derive their optimal safety-aware trajectory in a finite horizon.
The proposed control framework was validated by numerical simulations with multiple preceding HDVs showing that the generated control actions can ensure safe braking for the CAVs.
%future works
%Ongoing  efforts consider the complete solution that incorporates safety constraints associated with multiple preceding HDVs and lane changes. 
A direction for future research should focus on extending this framework to consider multi-lane traffic intersections with lane changing behavior of the HDVs.

%bibliography
\bibliography{references/ref_coordination, references/ref_mpc_acc, references/ids_bib_latest}       

\end{document}